\documentstyle[preprint]{elsarticle}
\begin{document}
\input{amssym}
\def\be{\begin{eqnarray}}
\def\ee{\end{eqnarray}}
\def\RR{{\Bbb R}}
\def\Span{\mathrm{Span}}
\title{Group analysis of three dimensional \\ Euler equations of gas dynamics}
\author[rvt]{Mehdi Nadjafikhah}
\address[rvt]{School of Mathematics, Iran University of Science and
Technology, \\ Narmak, Tehran, I.R.IRAN.} \ead{m
nadjafikhah@iust.ac.ir}
\date{}
\begin{abstract} In this paper, the equations governing the unsteady
flow of a perfect polytropic gas in three space dimensions are
considered. The basic similarity reductions for this system are
performed. Reduced equations and exact solutions associated with
the symmetries are obtained. This results is used to give the
invariance of system up to Galilean motions of space-times
$\RR^4$. Then, an optimal system of one-dimensional sub-algebras
for symmetry algebra of this system is given.
\end{abstract}
\maketitle
\section{Introduction}
The equations governing the unsteady flow of a perfect polytropic
gas in three space dimensions are
\be
&&\hspace{-1cm}  \frac{\partial u}{\partial t}+u\frac{\partial
u}{\partial x}+v\frac{\partial u}{\partial y}+w\frac{\partial
u}{\partial z}
+\frac{1}{q}\frac{\partial p}{\partial x}=0,\nonumber \\
&&\hspace{-1cm}  \frac{\partial v}{\partial t}+u\frac{\partial
v}{\partial x}+v\frac{\partial v}{\partial y}+w\frac{\partial
v}{\partial z}
+\frac{1}{q}\frac{\partial p}{\partial y}=0,\nonumber \\
&&\hspace{-1cm}  \frac{\partial w}{\partial t}+u\frac{\partial
w}{\partial x}+v\frac{\partial w}{\partial y}+w\frac{\partial
w}{\partial z}
+\frac{1}{q}\frac{\partial p}{\partial z}=0,\label{eq:1}\\
&&\hspace{-1cm}  \frac{\partial q}{\partial
t}+q\Big(\frac{\partial u}{\partial x}+\frac{\partial v}{\partial
y}+\frac{\partial w}{\partial z}\Big)+u\frac{\partial q}{\partial
x}+v\frac{\partial q}{\partial
y}+w\frac{\partial q}{\partial z}=0,\nonumber \\
&&\hspace{-1cm}  \frac{\partial p}{\partial t}+\gamma
p\Big(\frac{\partial u}{\partial x}+\frac{\partial v}{\partial
y}+\frac{\partial w}{\partial z}\Big)+u\frac{\partial p}{\partial
x}+v\frac{\partial p}{\partial y}+w\frac{\partial p}{\partial
z}=0.\nonumber
\ee
where $t$ is the time and $x$, $y$, $z$ are the space coordinates;
$q(t,x,y,z)$ is the density, $p(t,x,y,z)$ the pressure,
$u(t,x,y,z)$, $v(t,x,y,z)$ and $w(t,x,y,z)$ the velocity
components in the $x$, $y$ and $z$ directions, respectively, and
$\gamma\in\RR$ is the adiabatic index. These equations are called
"{\it three dimensional Euler equations of gas dynamics}"
\cite{SekSha} and section 6.3 of \cite{Ste}.

 In this paper, we consider the equations governing the
unsteady flow of a perfect polytropic gas in three space
dimensions. The basic similarity reductions for this system are
performed. Reduced equations and exact solutions associated with
the symmetries are obtained. We find an optimal system of
one-dimensional sub-algebras for symmetry algebra of this system.

In section 2 we find the full symmetry algebra of system
(\ref{eq:1}). Chapter 3 is devoted to the group-invariant
solutions to the system. The structure of full symmetry algebra
illustrated in section. The last section deals with a optimal
system of sub-algebras.
\section{Lie symmetry of the system}
The classical method for finding symmetry reductions of PDE is the
Lie group method of infinitesimal transformations. To apply the
classical method to (\ref{eq:1}), we consider the one-parameter
Lie group of infinitesimal transformations in $(t,x,y,z,u,v,p,q)$
given by
\be \tilde{t} &=& t+s.\xi_1(t,x,y,z,u,v,p,q)+O(s^2),\nonumber\\
\tilde{x} &=& x+s.\xi_2(t,x,y,z,u,v,p,q)+O(s^2),\nonumber\\
\tilde{y} &=& y+s.\xi_3(t,x,y,z,u,v,p,q)+O(s^2),\nonumber\\
\tilde{z} &=& z+s.\xi_4(t,x,y,z,u,v,p,q)+O(s^2),\nonumber\\
\tilde{u} &=& u+s.\eta_1(t,x,y,z,u,v,p,q)+O(s^2),\\
\tilde{v} &=& v+s.\eta_2(t,x,y,z,u,v,p,q)+O(s^2),\nonumber\\
\tilde{w} &=& w+s.\eta_3(t,x,y,z,u,v,p,q)+O(s^2),\nonumber\\
\tilde{p} &=& p+s.\eta_4(t,x,y,z,u,v,p,q)+O(s^2),\nonumber\\
\tilde{q} &=& q+s.\eta_5(t,x,y,z,u,v,p,q)+O(s^2). \nonumber\ee
where $s$ is the group parameter. One requires that this
transformation leaves the sub-manifold $S$ invariant. Let
$$X=\xi_1\,\partial_t+\xi_2\,\partial_x+\xi_3\,\partial_y+\xi_4\,\partial_z+\eta_1\,\partial_u+\eta_2\,\partial_v+\eta_3\,\partial_w+\eta_4\,\partial_p+\eta_5\,\partial_q$$
 be the corresponding infinitesimal transformation; where
coefficients are real-valued $C^\infty-$functions of
$t,x,y,z,u,v,w,p,q$. The infinitesimal criterion of invariance
(\ref{eq:1}) will thus involve $t,x,y,z,u,v,w,p,q$ and the
derivatives of $u,v,w,p,q$ with respect to $t,x,y,z$, as well as
$\xi_i$, $i=1,2,3,4$, and $\eta_j$, $j=1,2,3,4,5$, and their
partial derivatives. After eliminating any dependencies among the
derivatives of the m's caused by the system itself, we can then
equate the coefficients of the remaining unconstrained partial
derivatives of u to zero. This will result in a large number of
elementary partial differential equations for the coefficient
functions $\xi_i$, $\eta_j$ of the infinitesimal generator, called
the determining equations for the symmetry group of the given
system. In this case, we find the determining equations for the
symmetry group of the system (\ref{eq:1}) to be the following:
\be
&&\hspace{-1cm} \xi_{1,t^2}=\xi_{1,x}=\xi_{1,y}=\xi_{1,z}=\xi_{1,u}=\xi_{1,v}=\xi_{1,w}=\xi_{1,q}=\xi_{1,p}=0,\nonumber\\
&&\hspace{-1cm} \xi_{2,x}=\xi_{2,u}=\xi_{2,v}=\xi_{2,w}=\xi_{2,q}=\xi_{2,p}=\xi_{2,t^2}=\xi_{2,ty}=\xi_{2,y^2}=0,\nonumber\\
&&\hspace{-1cm} \xi_{3,y}=\xi_{3,u}=\xi_{3,v}=\xi_{3,w}=\xi_{3,q}=\xi_{3,p}=\xi_{3,t^2}=0,\label{eq:3}\\
&&\hspace{-1cm} \xi_{4,z}=\xi_{4,w}=\xi_{4,q}=\xi_{4,p}=0,\nonumber\\
&&\hspace{-1cm} \xi_{4,t^2}=\xi_{4,tx}=\xi_{4,ty}=\xi_{4,x^2}=\xi_{4,xy}=\xi_{4,y^2}=0,\nonumber\\
&&\hspace{-1cm}
\xi_{3,x}=\xi_{2,y},\,\xi_{3,z}+\xi_{4,y}=0,\,\xi_{4,v}=\xi_{4,u},\,\xi_{2,z}+\xi_{4,x}=0,\nonumber
\ee
and
\be
&&\hspace{-1cm} \eta_1=-u\xi_{1,t}+\xi_{2,t}+v\xi_{2,y}-w\xi_{4,x},\nonumber\\
&&\hspace{-1cm} \eta_2=-w\xi_{4,y}+\xi_{3,t}-u\xi_{2,y}-v\xi_{1,t},\nonumber\\
&&\hspace{-1cm} \eta_3=u\xi_{4,x}+\xi_{4,t}+v\xi_{4,y}-w\xi_{1,t},\label{eq:4}\\
&&\hspace{-1cm} \eta_4=q\eta_{5}/p+2q\xi_{1,t},\nonumber\\
&&\hspace{-1cm} \eta_{5,t}=\eta_{5,x}=\eta_{5,y}=\eta_{5,z}=\eta_{5,u}=\eta_{5,v}=\eta_{5,w}=\eta_{5,q}=0,\nonumber\\
&&\hspace{-1cm} \eta_{5}=p\eta_{5,p}.\nonumber \ee
First, equations (\ref{eq:3}) require that $\xi_i$s are just
\be
&&\hspace{-1cm} \xi_1=a_1t+a_2,\nonumber \\
&&\hspace{-1cm} \xi_2=a_3t+a_4x+a_5y+a_6z+a_7,\nonumber \\
&&\hspace{-1cm} \xi_3=a_8t-a_5x+a_4y+a_9z+a_{10},\label{eq:5}\\
&&\hspace{-1cm} \xi_4=a_{11}t-a_6x-a_9y+a_4z+a_{12},\nonumber
\ee
where $a_1,\cdots,a_{12}$ are arbitrary constants. Then, equations
(\ref{eq:4}) and (\ref{eq:5}) requires that
\be
&&\hspace{-1cm} \eta_1=(a_4-a_1)u+a_5v-a_6w+a_3,\nonumber \\
&&\hspace{-1cm} \eta_2=-a_5u+(a_4-a_1)v+a_4w+a_8,\nonumber \\
&&\hspace{-1cm} \eta_3=-a_6u-a_9v+(a_4-a_1)w+a_{11},\label{eq:6}\\
&&\hspace{-1cm} \eta_4=(2a_1-2a_4+a_{13})q,\nonumber \\
&&\hspace{-1cm} \eta_5=a_{13}p,\nonumber
\ee
where $a_{13}$ is an arbitrary constant. Therefore
\paragraph{Theorem 1.}
{\em The Lie algebra of infinitesimal symmetries of the system
(\ref{eq:1}) is spanned by the 13 vector fields
\be
\begin{array}{lll}
X_1 = \partial_x,&&\\
X_2 = \partial_y,&&\hspace{-3cm}\mbox{(shifts of origin)}\\
X_3 = \partial_z,&&\\
X_4 = \partial_t,&&\\ \\
X_5 = t\,\partial_x+\partial_u,&&\\
X_6 = t\,\partial_y+\partial_v,&&\hspace{-3cm}\mbox{(uniformly frame motions)}\\
X_7 = t\,\partial_z+\partial_w,&&\\ \\
X_8 = t\,\partial_t-u\,\partial_u-v\,\partial_v-w\,\partial_w+2q\,\partial_q,\\
X_9 =
x\,\partial_x+y\,\partial_y+z\,\partial_z+u\,\partial_u+v\,\partial_v+w\,\partial_w-2r\,\partial_r,\label{eq:7}\\
\\
X_{10} = y\,\partial_x-x\,\partial_y+v\,\partial_u-u\,\partial_v,&&\\
X_{11} = -z\,\partial_y+y\,\partial_z-w\,\partial_v+v\,\partial_w,&&\hspace{-3cm}\mbox{(rotations of reference frame)}\\
X_{12} =
-z\,\partial_x+x\,\partial_z-w\,\partial_u+u\,\partial_w,&&\\ \\
X_{13} = q\,\partial_q+p\,\partial_p,
\end{array}
\ee
These infinitesimal symmetries will generate a Lie algebra $\goth
g$ over the field of real or complex numbers. The commutator table
of Lie algebra $\goth g$ for (\ref{eq:1}) is given below, where
the entry in the $i$th row and $j$th column is defined as
$[X_i,X_j]=X_i.X_j-X_j.X_i$, $i,j=1,\cdots,13$.
\begin{eqnarray*}
&& \hspace{-0.7cm}\small \begin{array}{|c|ccccccccccccc|}
\hline\   \!\!&\!\!X_1\!\!&\!\!X_2\!\!&\!\!X_3\!\!&\!\!X_4\!\!&\!\!X_5\!\!&\!\!X_6\!\!&\!\!X_7\!\!&\!\!X_8\!\!&\!\!X_9\!\!&\!\!X_{10}\!\!&\!\!X_{11}\!\!&\!\!X_{12}\!\!&\!\!X_{13}\\
\hline\ X_1\!\!&\!\!0\!\!&\!\!0\!\!&\!\!0\!\!&\!\!0\!\!&\!\!0\!\!&\!\!0\!\!&\!\!0\!\!&\!\!0\!\!&\!\!X_1\!\!&\!\!-X_2\!\!&\!\!0\!\!&\!\!X_3\!\!&\!\!0\\
X_2\!\!&\!\!0\!\!&\!\!0\!\!&\!\!0\!\!&\!\!0\!\!&\!\!0\!\!&\!\!0\!\!&\!\!0\!\!&\!\!0\!\!&\!\!X_2\!\!&\!\!X_1\!\!&\!\!X_3\!\!&\!\!0\!\!&\!\!0\\
X_3\!\!&\!\!0\!\!&\!\!0\!\!&\!\!0\!\!&\!\!0\!\!&\!\!0\!\!&\!\!0\!\!&\!\!0\!\!&\!\!0\!\!&\!\!X_3\!\!&\!\!0\!\!&\!\!-X_2\!\!&\!\!-X_1\!\!&\!\!0\\
X_4\!\!&\!\!0\!\!&\!\!0\!\!&\!\!0\!\!&\!\!0\!\!&\!\!X_1\!\!&\!\!X_2\!\!&\!\!X_3\!\!&\!\!X_4\!\!&\!\!0\!\!&\!\!0\!\!&\!\!0\!\!&\!\!0\!\!&\!\!0\\
X_5\!\!&\!\!0\!\!&\!\!0\!\!&\!\!0\!\!&\!\!-X_1\!\!&\!\!0\!\!&\!\!0\!\!&\!\!0\!\!&\!\!-X_5\!\!&\!\!X_5\!\!&\!\!-X_6\!\!&\!\!0\!\!&\!\!X_7\!\!&\!\!0\\
X_6\!\!&\!\!0\!\!&\!\!0\!\!&\!\!0\!\!&\!\!-X_2\!\!&\!\!0\!\!&\!\!0\!\!&\!\!0\!\!&\!\!-X_6\!\!&\!\!X_6\!\!&\!\!X_5\!\!&\!\!X_7\!\!&\!\!0\!\!&\!\!0\\
X_7\!\!&\!\!0\!\!&\!\!0\!\!&\!\!0\!\!&\!\!-X_3\!\!&\!\!0\!\!&\!\!0\!\!&\!\!0\!\!&\!\!-X_7\!\!&\!\!X_7\!\!&\!\!0\!\!&\!\!-X_6\!\!&\!\!-X_5\!\!&\!\!0\\
X_8\!\!&\!\!0\!\!&\!\!0\!\!&\!\!0\!\!&\!\!-X_4\!\!&\!\!X_5\!\!&\!\!X_6\!\!&\!\!X_7\!\!&\!\!0\!\!&\!\!0\!\!&\!\!0\!\!&\!\!0\!\!&\!\!0\!\!&\!\!0\\
X_9\!\!&\!\!-X_1\!\!&\!\!-X_2\!\!&\!\!-X_3\!\!&\!\!0\!\!&\!\!-X_5\!\!&\!\!-X_6\!\!&\!\!-X_7\!\!&\!\!0\!\!&\!\!0\!\!&\!\!0\!\!&\!\!0\!\!&\!\!0\!\!&\!\!0\\
X_{10}\!\!&\!\!X_2\!\!&\!\!-X_1\!\!&\!\!0\!\!&\!\!0\!\!&\!\!X_6\!\!&\!\!-X_5\!\!&\!\!0\!\!&\!\!0\!\!&\!\!0\!\!&\!\!0\!\!&\!\!-X_{12}\!\!&\!\!X_{11}\!\!&\!\!0\\
X_{11}\!\!&\!\!0\!\!&\!\!-X_3\!\!&\!\!X_2\!\!&\!\!0\!\!&\!\!0\!\!&\!\!-X_7\!\!&\!\!X_6\!\!&\!\!0\!\!&\!\!0\!\!&\!\!X_{12}\!\!&\!\!0\!\!&\!\!-X_{10}\!\!&\!\!0\\
X_{12}\!\!&\!\!-X_3\!\!&\!\!0\!\!&\!\!X_1\!\!&\!\!0\!\!&\!\!-X_7\!\!&\!\!0\!\!&\!\!X_5\!\!&\!\!0\!\!&\!\!0\!\!&\!\!-X_{11}\!\!&\!\!X_{10}\!\!&\!\!0\!\!&\!\!0\\
X_{13}\!\!&\!\!0\!\!&\!\!0\!\!&\!\!0\!\!&\!\!0\!\!&\!\!0\!\!&\!\!0\!\!&\!\!0\!\!&\!\!0\!\!&\!\!0\!\!&\!\!0\!\!&\!\!0\!\!&\!\!0\!\!&\!\!0\\
\hline
\end{array}
\end{eqnarray*}}
\section{Reduction of the system}
Now, we look at the group-invariant solutions to the system
(\ref{eq:1}).

\medskip The system (\ref{eq:1}) is a sub-manifold $M$ of jet space
$J^2(\RR^4,\RR^5)$ defined by the following equations
\be
&&\hspace{-1cm} u_t+uu_x+vu_y+wu_z +p_x/q=0,\nonumber \\
&&\hspace{-1cm} v_t+uv_x+vv_y+wv_z+p_y/q=0,\nonumber \\
&&\hspace{-1cm} w_t+uw_x+vw_y+ww_z+p_z/q=0,\label{eq:2}\\
&&\hspace{-1cm} q_t+q \left(u_x+v_y+w_z\right)+uq_x+vq_y+wq_z=0,\nonumber \\
&&\hspace{-1cm}
p_t+\gamma\,p\left(u_x+v_y+w_z\right)+up_x+vp_y+wp_z=0.\nonumber
\ee

Doing as section 3.1 of \cite{Olv2} and find (in a sense) the most
general group-invariant solutions to the system (\ref{eq:1}).
\paragraph{Theorem 2.}
{\em The one-parameter groups $g_i(s):M\to M$ generated by the
$X_i$, $i=1,\cdots,13$ are given in the following table:
\be
g_1 &:& \big(t,x+s,y,z,u,v,w,p,q\big)\nonumber\\
g_2 &:& \big(t,x,y+s,z,u,v,w,p,q\big)\nonumber\\
g_3 &:& \big(t,x,y,z+s,u,v,w,p,q\big)\nonumber\\
g_4 &:& \big(t+s,x,y,z,u,v,w,p,q\big)\nonumber\\
g_5 &:& \big(t,x+st,y,z,u+s,v,w,p,q\big) \nonumber\\
g_6 &:& \big(t,x,y+st,z,u,v+s,w,p,q\big) \nonumber\\
g_7 &:& \big(t,x,y,z+st,u,v,w+s,p,q\big) \label{eq:8}\\
g_8 &:& \big(e^st,x,y,z,e^{-s}u,e^{-s}v,e^{-s}w,p,e^{2s}q\big) \nonumber\\
g_9 &:& \big(t,e^sx,e^sy,e^sz,e^su,e^sv,e^sw,p,e^{-2s}q\big) \nonumber\\
g_{10} &:& \big(t,x\cos s+y\sin s,y\cos s-x\sin s,z,\nonumber\\
&& \hspace{3cm} v\sin s+u\cos s,v\cos s-u\sin s,w,p,q\big) \nonumber\\
g_{11} &:& \big(t,x,y\cos s-z\sin s,z\cos s+y\sin s,\nonumber\\
&& \hspace{3cm} u,v\cos s-w\sin s,w\cos s+v\sin s,p,q\big) \nonumber\\
g_{12} &:& \big(t,x\cos s-z\sin s,y,x\sin s+z\cos s,\nonumber\\
&& \hspace{3cm} u\cos s-w\sin s,v,w\cos s+u\sin s,p,q\big) \nonumber\\
g_{13} &:& \big(t,x,y,z,u,v,w,e^sp,e^sq\big) \nonumber
\ee
Where, entries give the transformed point
$$\exp(sX_i)(t,x,y,z,u,v,w,p,q) =
(\tilde{t},\tilde{x},\tilde{y},\tilde{z},\tilde{u},\tilde{v},\tilde{w},\tilde{p},\tilde{q}).$$}

\medskip Since each group $g_i$, $i=1,\cdots,13$, is a symmetry
group of (\ref{eq:1}), then
\paragraph{Theorem 3.}
{\em Let $i=1,\cdots,13$, $s\in\RR$, and $\lambda\in\RR-\{0\}$. If
\begin{eqnarray*}
&& u=U(t,x,y,z),\;v=V(t,x,y,z),\;w=W(t,x,y,z),\\
&& p=P(t,x,y,z),\;q=Q(t,x,y,z),
\end{eqnarray*}
is a solution of the system (\ref{eq:1}), so are the functions
\begin{eqnarray*}
&& u^{(i)}=U^{(i)}(t,x,y,z),\;v^{(i)}=V^{(i)}(t,x,y,z),\;w^{(i)}=W^{(i)}(t,x,y,z),\\
&& p^{(i)}=P^{(i)}(t,x,y,z),\;q^{(i)}=Q^{(i)}(t,x,y,z),
\end{eqnarray*}
where
\be
\begin{array}{llll}
u^{(1)}\,=\,U\big(t,x+s,y,z\big), &&& u^{(2)}\,=\,U\big(t,x,y+s,z\big),\\
v^{(1)}\,=\,V\big(t,x+s,y,z\big), &&& v^{(2)}\,=\,V\big(t,x,y+s,z\big),\\
w^{(1)}\,=\,W\big(t,x+s,y,z\big), &&& w^{(2)}\,=\,W\big(t,x,y+s,z\big),\\
p^{(1)}\,=\,P\big(t,x+s,y,z\big), &&& p^{(2)}\,=\,P\big(t,x,y+s,z\big),\\
q^{(1)}\,=\,Q\big(t,x+s,y,z\big). &&&
q^{(2)}\,=\,Q\big(t,x,y+s,z\big).
\end{array} \nonumber
\ee
\be
\begin{array}{llll}
u^{(3)}\,=\,U\big(t,x,y,z+s\big), &&& u^{(4)}\,=\,U\big(t+s,x,y,z\big),\\
v^{(3)}\,=\,V\big(t,x,y,z+s\big), &&& v^{(4)}\,=\,V\big(t+s,x,y,z\big),\\
w^{(3)}\,=\,W\big(t,x,y,z+s\big), &&& w^{(4)}\,=\,W\big(t+s,x,y,z\big),\\
p^{(3)}\,=\,P\big(t,x,y,z+s\big), &&& p^{(4)}\,=\,P\big(t+s,x,y,z\big),\\
q^{(3)}\,=\,Q\big(t,x,y,z+s\big). &&&
q^{(4)}\,=\,Q\big(t+s,x,y,z\big).
\end{array} \nonumber
\ee
\be
\begin{array}{llll}
u^{(5)}\,=\,U\big(t,x+st,y,z\big)-s, &&& u^{(6)}\,=\,U\big(t,x,y+st,z\big),\\
v^{(5)}\,=\,V\big(t,x+st,y,z\big),   &&& v^{(6)}\,=\,V\big(t,x,y+st,z\big)-s,\\
w^{(5)}\,=\,W\big(t,x+st,y,z\big),   &&& w^{(6)}\,=\,W\big(t,x,y+st,z\big),\\
p^{(5)}\,=\,P\big(t,x+st,y,z\big),   &&& p^{(6)}\,=\,P\big(t,x,y+st,z\big),\\
q^{(5)}\,=\,Q\big(t,x+st,y,z\big).   &&&
q^{(6)}\,=\,Q\big(t,x,y+st,z\big).
\end{array} \nonumber
\ee
\be
\begin{array}{llll}
u^{(7)}\,=\,U\big(t,x,y,z+st\big),   &&& u^{(8)}\,=\,\lambda.U\big(\lambda.t,x,y,z\big),\\
v^{(7)}\,=\,V\big(t,x,y,z+st\big),   &&& v^{(8)}\,=\,\lambda.V\big(\lambda.t,x,y,z\big),\\
w^{(7)}\,=\,W\big(t,x,y,z+st\big)-s, &&& w^{(8)}\,=\,\lambda.W\big(\lambda.t,x,y,z\big),\\
p^{(7)}\,=\,P\big(t,x,y,z+st\big),   &&& p^{(8)}\,=\,P\big(\lambda.t,x,y,z\big),\\
q^{(7)}\,=\,Q\big(t,x,y,z+st\big).   &&&
q^{(8)}\,=\,\frac{1}{\lambda^2}.Q\big(\lambda.t,x,y,z\big).
\end{array} \nonumber
\ee
\be
\begin{array}{llll}
u^{(9)}\,=\,\frac{1}{\lambda}.U\big(t,\lambda.x,\lambda.y,\lambda.z\big), &&& u^{(13)}\,=\,U\big(t,x,y,z\big),\\
v^{(9)}\,=\,\frac{1}{\lambda}.V\big(t,\lambda.x,\lambda.y,\lambda.z\big), &&& v^{(13)}\,=\,V\big(t,x,y,z\big),\\
w^{(9)}\,=\,\frac{1}{\lambda}.W\big(t,\lambda.x,\lambda.y,\lambda.z\big), &&& w^{(13)}\,=\,W\big(t,x,y,z\big),\\
p^{(9)}\,=\,P\big(t,\lambda.x,\lambda.y,\lambda.z\big),                   &&& p^{(13)}\,=\,\lambda.P\big(t,x,y,z\big),\\
q^{(9)}\,=\,\lambda^2.Q\big(t,\lambda.x,\lambda.y,\lambda.z\big).
&&& q^{(13)}\,=\,\lambda.Q\big(t,x,y,z\big).
\end{array} \nonumber
\ee
\be
u^{(10)}&=&\cos s.U\big(t,x\cos s+y\sin s,y\cos s-x\sin s,z\big)\nonumber\\
&& -\sin s.V\big(t,x\cos s+y\sin s,y\cos s-x\sin s,z\big),\nonumber\\
v^{(10)}&=&\sin s.U\big(t,x\cos s+y\sin s,y\cos s-x\sin s,z\big)\nonumber\\
&&+\cos s.V\big(t,x\cos s+y\sin s,y\cos s-x\sin s,z\big),\nonumber\\
w^{(10)}&=&W\big(t,x\cos s+y\sin s,y\cos s-x\sin s,z\big),\nonumber\\
p^{(10)}&=&P\big(t,x\cos s+y\sin s,y\cos s-x\sin s,z\big),\nonumber\\
q^{(10)}&=&Q\big(t,x\cos s+y\sin s,y\cos s-x\sin
s,z\big).\nonumber \ee
\be
u^{(11)}&=&U\big(t,x,y\cos s-z\sin s,z\cos s+y\sin s\big),\nonumber\\
v^{(11)}&=&\cos s.V\big(t,x,y\cos s-z\sin s,z\cos s+y\sin s\big)\nonumber\\
&&+\sin s.W\big(t,x,y\cos s-z\sin s,z\cos s+y\sin s\big),\nonumber\\
w^{(11)}&=&\cos s.W\big(t,x,y\cos s-z\sin s,z\cos s+y\sin s\big)\nonumber\\
&&-\sin s.V\big(t,x,y\cos s-z\sin s,z\cos s+y\sin
s\big),\nonumber\\
p^{(11)}&=&P\big(t,x,y\cos s-z\sin s,z\cos s+y\sin s\big),\nonumber\\
q^{(11)}&=&Q\big(t,x,y\cos s-z\sin s,z\cos s+y\sin
s\big).\nonumber \ee
\be
u^{(12)}&=&\cos s.U\big(t,x\cos s-z\sin s,y,x\sin s+z\cos s\big)\nonumber\\
&&+\sin s.W\big(t,x\cos s-z\sin s,y,x\sin s+z\cos s\big),\nonumber\\
v^{(12)}&=&V\big(t,x\cos s-z\sin s,y,x\sin s+z\cos s\big),\nonumber\\
w^{(12)}&=&\cos s.W\big(t,x\cos s-z\sin s,y,x\sin s+z\cos s\big)\nonumber\\
&&-\sin s.U\big(t,x\cos s-z\sin s,y,x\sin s+z\cos
s\big),\nonumber\\
p^{(12)}&=&P\big(t,x\cos s-z\sin s,y,x\sin s+z\cos s\big),\nonumber\\
q^{(12)}&=&Q\big(t,x\cos s-z\sin s,y,x\sin s+z\cos
s\big)\nonumber. \ee}
\section{Structure of Lie algebra $\goth g$}
In this section, we determine the structure of full symmetry
algebra $\goth g$. of system (\ref{eq:1}).

The center $\goth z$ of $\goth g$ is $\Span_\RR\{X_{13}\}$.
Therefore, the quotient algebra ${\goth g}_1:={\goth g}/{\goth z}$
is $\Span_\RR\{Y_1,\cdots,Y_{12}\}$; where $Y_i:=X_i+{\goth z}$,
$i=1,\cdots,12$. The commutator table of Lie algebra ${\goth g}_1$
is given below, where the entry in the $i$th row and $j$th column
is defined as $[Y_i,Y_j]=Y_i.Y_j-Y_j.Y_i$, $i,j=1,\cdots,12$.
\begin{eqnarray*}
&& \hspace{-1.3cm}\small \begin{array}{|c|cccccccccccc|}
\hline\   \!\!&\!\!Y_1\!\!&\!\!Y_2\!\!&\!\!Y_3\!\!&\!\!Y_4\!\!&\!\!Y_5\!\!&\!\!Y_6\!\!&\!\!Y_7\!\!&\!\!Y_8\!\!&\!\!Y_9\!\!&\!\!Y_{10}\!\!&\!\!Y_{11}\!\!&\!\!Y_{12}\\
\hline\ Y_1\!\!&\!\!0\!\!&\!\!0\!\!&\!\!0\!\!&\!\!0\!\!&\!\!0\!\!&\!\!0\!\!&\!\!0\!\!&\!\!0\!\!&\!\!Y_1\!\!&\!\!-Y_2\!\!&\!\!0\!\!&\!\!Y_3\\
Y_2\!\!&\!\!0\!\!&\!\!0\!\!&\!\!0\!\!&\!\!0\!\!&\!\!0\!\!&\!\!0\!\!&\!\!0\!\!&\!\!0\!\!&\!\!Y_2\!\!&\!\!Y_1\!\!&\!\!Y_3\!\!&\!\!0\\
Y_3\!\!&\!\!0\!\!&\!\!0\!\!&\!\!0\!\!&\!\!0\!\!&\!\!0\!\!&\!\!0\!\!&\!\!0\!\!&\!\!0\!\!&\!\!Y_3\!\!&\!\!0\!\!&\!\!-Y_2\!\!&\!\!-Y_1\\
Y_4\!\!&\!\!0\!\!&\!\!0\!\!&\!\!0\!\!&\!\!0\!\!&\!\!Y_1\!\!&\!\!Y_2\!\!&\!\!Y_3\!\!&\!\!Y_4\!\!&\!\!0\!\!&\!\!0\!\!&\!\!0\!\!&\!\!0\\
Y_5\!\!&\!\!0\!\!&\!\!0\!\!&\!\!0\!\!&\!\!-Y_1\!\!&\!\!0\!\!&\!\!0\!\!&\!\!0\!\!&\!\!-Y_5\!\!&\!\!Y_5\!\!&\!\!-Y_6\!\!&\!\!0\!\!&\!\!Y_7\\
Y_6\!\!&\!\!0\!\!&\!\!0\!\!&\!\!0\!\!&\!\!-Y_2\!\!&\!\!0\!\!&\!\!0\!\!&\!\!0\!\!&\!\!-Y_6\!\!&\!\!Y_6\!\!&\!\!Y_5\!\!&\!\!Y_7\!\!&\!\!0\\
Y_7\!\!&\!\!0\!\!&\!\!0\!\!&\!\!0\!\!&\!\!-Y_3\!\!&\!\!0\!\!&\!\!0\!\!&\!\!0\!\!&\!\!-Y_7\!\!&\!\!Y_7\!\!&\!\!0\!\!&\!\!-Y_6\!\!&\!\!-Y_5\\
Y_8\!\!&\!\!0\!\!&\!\!0\!\!&\!\!0\!\!&\!\!-Y_4\!\!&\!\!Y_5\!\!&\!\!Y_6\!\!&\!\!Y_7\!\!&\!\!0\!\!&\!\!0\!\!&\!\!0\!\!&\!\!0\!\!&\!\!0\\
Y_9\!\!&\!\!-Y_1\!\!&\!\!-Y_2\!\!&\!\!-Y_3\!\!&\!\!0\!\!&\!\!-Y_5\!\!&\!\!-Y_6\!\!&\!\!-Y_7\!\!&\!\!0\!\!&\!\!0\!\!&\!\!0\!\!&\!\!0\!\!&\!\!0\\
Y_{10}\!\!&\!\!Y_2\!\!&\!\!-Y_1\!\!&\!\!0\!\!&\!\!0\!\!&\!\!Y_6\!\!&\!\!-Y_5\!\!&\!\!0\!\!&\!\!0\!\!&\!\!0\!\!&\!\!0\!\!&\!\!-Y_{12}\!\!&\!\!Y_{11}\\
Y_{11}\!\!&\!\!0\!\!&\!\!-Y_3\!\!&\!\!Y_2\!\!&\!\!0\!\!&\!\!0\!\!&\!\!-Y_7\!\!&\!\!Y_6\!\!&\!\!0\!\!&\!\!0\!\!&\!\!Y_{12}\!\!&\!\!0\!\!&\!\!-Y_{10}\\
Y_{12}\!\!&\!\!-Y_3\!\!&\!\!0\!\!&\!\!Y_1\!\!&\!\!0\!\!&\!\!-Y_7\!\!&\!\!0\!\!&\!\!Y_5\!\!&\!\!0\!\!&\!\!0\!\!&\!\!-Y_{11}\!\!&\!\!Y_{10}\!\!&\!\!0\\
\hline
\end{array}
\end{eqnarray*}
The algebra ${\goth g}_1$ is not solvable, because
\[ {\goth g}_1^{(1)}:=[{\goth g}_1,{\goth
g}_1]=\Span_\RR\{Y_1,\cdots,Y_{10}\},\;\;\mbox{and}\;\; {\goth
g}_1^{(2)}:=[{\goth g}_1^{(1)},{\goth g}_1^{(1)}]={\goth
g}_1^{(1)}.
\]
Lie algebra ${\goth g}_1$ admits a Levi-decomposition into the
semi-direct sum $${\goth g}_1={\goth r}\ltimes{\goth s},$$ where
${\goth r}=\Span_\RR\{Y_1,\cdots,Y_9\}$ is the radical of ${\goth
g}_1$ (the largest solvable ideal contained in ${\goth g}_1$) and
${\goth s}=\Span_\RR\{Y_{10},Y_{11},Y_{12}\}$ is a semi-simple
subalgebra of ${\goth g}_1$. The radical $\goth r$ is uniquely
defined but the semi-simple subalgebra $\goth s$ is not. $\goth s$
as an algebra is isomorphic to $3-$dimensional special orthogonal
algebra ${\goth so}(3):=\{A\in\mathrm{Mat}(3\times
3;\RR)\,|\,A^t=-A\}$, which is a simple algebra.

The radical ${\goth r}$ is solvable, with and
$${\goth r}^{(1)}=\Span_\RR\{Y_1,\cdots,Y_7\},\;\;\;\; {\goth
r}^{(2)}=\Span_\RR\{Y_1,Y_2,Y_3\},\;\;\;\; {\goth
r}^{(3)}=\{0\}.$$ It is semi-direct sum ${\goth r}={\goth
r}^{(1)}\ltimes {\goth g}_2$, where ${\goth
g}_2=\Span_\RR\{Y_8,Y_9\}$ is isomorphic to the two dimensional
Abelian Lie algebra $\RR^2$; i.e. $${\goth r}={\goth
r}^{(1)}\ltimes \RR^2.$$

${\goth r}^{(1)}$ is semi-direct sum ${\goth r}^{(1)}={\goth
r}^{(2)}\ltimes {\goth g}_3,$ where ${\goth
g}_3=\Span_\RR\{Y_4,\cdots,Y_7\}$ is isomorphic to the four
dimensional Abelian Lie algebra $\RR^4$, and ${\goth r}^{(2)}$ is
isomorphic to the three dimensional Abelian Lie algebra $\RR^3$;
i.e.
$${\goth r}^{(1)}=\RR^3\ltimes \RR^4.$$
\paragraph{Theorem 4.}
{\em The symmetry algebra $\goth g$ of system (\ref{eq:1}) have
the following structure $${\goth
g}\cong\RR\times(((\RR^3\ltimes\RR^4)\ltimes\RR^2)\ltimes\goth{
so}(3)).$$}
\paragraph{Theorem 5.}
{\em There is a normal Lie-subalgebra of symmetry Lie-group of
system (\ref{eq:1}) which is isomorphic to the Galilean group of
$\RR^4$: $${\rm Gal}(4,{\Bbb R})=\left\{ \left(\begin{array}{ccc} 1 & {\bf 0} & s \\ {\bf v} & {\bf R} & {\bf y} \\
0 & {\bf 0} & 1
\end{array}\right)\;\Big|\; {\bf R}\in{\rm
O}(3,{\Bbb R}),\; s\in{\Bbb R},\; \mbox{and}\;\; {\bf y},{\bf
v}\in{\Bbb R}^3 \;\right\}.$$ Therefore, the system (\ref{eq:1})
is invariant up to Galilean motions of space-times $\RR^4$ (see
\cite{NadMeh}).}

\medskip \noindent {\it Proof.} Let
${\goth b}=\Span_\RR\{X_1,\cdots,X_7,X_{10},X_{11},X_{12}\}$. Lie
algebra structure of ${\goth b}$ is as $\goth{gal}(3)$, the lie
algebra of Galilean group ${\rm Gal}(4,{\Bbb R})$, and $\goth b$
an ideal of $\goth g$. Therefore, there is a Lie-subgroup of $G$
such that its Lie-algebra is $\goth b$, by the Theorem 2.53 of
Olver \cite{Olv2}. \hfill\ $\Box$
\paragraph{Conclusion 1.}
{\em $\goth g$ is semi-direct sum of $$\goth{
gal}(4)\cong\Span_\RR\{X_1,\cdots,X_7,X_{10},X_{11},X_{12}\}$$ and
$3-$dimensional Abelian algebra
$\RR^3\cong\Span\{X_8,X_9,X_{13}\}$.}
\section{Optimal system of sub-algebras}
As is well known, the Lie group theoretic method plays an
important role in finding exact solutions and performing symmetry
reductions of differential equations. Since any linear combination
of infinitesimal generators is also an infinitesimal generator,
there are always infinitely many different symmetry subgroups for
the differential equation. So, a mean of determining which
subgroups would give essentially different types of solutions is
necessary and significant for a complete understanding of the
invariant solutions. As any transformation in the full symmetry
group maps a solution to another solution, it is sufficient to
find invariant solutions which are not related by transformations
in the full symmetry group, this has led to the concept of an
optimal system \cite{Ovs}. The problem of finding an optimal
system of subgroups is equivalent to that of finding an optimal
system of subalgebras. For one-dimensional subalgebras, this
classification problem is essentially the same as the problem of
classifying the orbits of the adjoint representation. This problem
is attacked by the naive approach of taking a general element in
the Lie algebra and subjecting it to various adjoint
transformations so as to simplify it as much as possible. The idea
of using the adjoint representation to classify group-invariant
solutions was due to \cite{Ovs} and \cite{Olv1}.

The adjoint action is given by the Lie series
\be \mathrm{Ad}(\exp(sX_i)X_j) =
X_j-s[X_i,X_j]+\frac{s^2}{2}[X_i,[X_i,X_j]]-\cdots,\ee
where $[X_i,X_j]$ is the commutator for the Lie algebra, $s$ is a
parameter, and $i,j=1,\cdots,13$. We can write the adjoint action
for the Lie algebra $\goth g$, and show that
\paragraph{Theorem 6.}
{\em A one-dimensional optimal system of (\ref{eq:1}) is given by
\begin{itemize}
\item[1)$\;\;$] $X=X_8+\sum_{i=9}^{13}a_iX_i$,
\item[2)$\;\;$] $X=a_7X_7+X_8+X_9+ \sum_{i=10}^{13}a_iX_i$,
\item[3)$\;\;$] $X=a_6X_6+a_7X_7+X_8+X_9+ \sum_{i=11}^{13}a_iX_i$,
\item[4)$\;\;$] $X=a_5X_5+X_8+X_9+a_{11}X_i+a_{13}X_{13}$,
\item[5)$\;\;$] $X=\sum_{i=5}^7a_iX_i+X_8+X_9+a_{13}X_{13}$,
\item[6)$\;\;$] $X=a_4X_4+X_9+\sum_{i=10}^{13}a_iX_i$,
\item[7)$\;\;$] $X=X_4+a_5X_5+a_8X_8+\sum_{i=10}^{13}a_iX_i$,
\item[8)$\;\;$] $X=X_4+a_6X_6+a_8X_8+a_{10}X_{10}+a_{12}X_{12}+a_{13}X_{13}$,
\item[9)$\;\;$] $X=X_4+a_7X_7+a_8X_8+a_{10}X_{10}+a_{12}X_{12}+a_{13}X_{13}$,
\item[10)] $X=X_4+\sum_{i=5}^8a_iX_i+a_{13}X_{13}$,
\item[11)] $X=X_4+a_7X_7+\sum_{i=10}^{13}a_iX_i$,
\item[12)] $X=X_7+\sum_{i=10}^{13}a_iX_i$,
\item[13)] $X=a_3X_3+\sum_{i=10}^{13}a_iX_i$,
\item[14)] $X=X_4+a_6X_6+\sum_{i=11}^{13}a_iX_i$,
\item[15)] $X=X_4+a_5X_5+a_{11}X_{11}+a_{13}X_{13}$,
\item[16)] $X=X_4+\sum_{i=5}^7a_iX_i+a_{13}X_{13}$,
\item[17)] $X=X_6+\sum_{i=11}^{13}a_iX_i$,
\item[18)] $X=a_2X_2+\sum_{i=11}^{13}a_iX_i$,
\item[19)] $X=a_2X_2+a_{12}X_{12}+a_{13}X_{13}$,
\item[20)] $X=X_5+a_{11}X_{11}+a_{13}X_{13}$,
\item[21)] $X=a_2X_2+a_3X_3+X_5+a_{13}X_{13}$,
\item[22)] $X=a_1X_1+X_{11}+a_{13}X_{13}$,
\item[23)] $X=a_1X_1+a_3X_3+X_6+a_{13}X_{13}$,
\item[24)] $X=a_1X_1+a_2X_2+a_3X_3+a_{13}X_{13}$,
\item[25)] $X=X_5+a_{10}X_{10}+a_{11}X_{11}+a_{13}X_{13}$,
\item[26)] $X=a_3X_3+X_{10}+a_{11}X_{11}+a_{13}X_{13}$.
\end{itemize}}

\medskip \noindent {\it Proof:} $F(sX_i):{\goth g}\to{\goth g}$ defined by
$X\mapsto\mathrm{Ad}(\exp(sX_i)X)$ is a linear map, for
$i=1,\cdots,13$. The matrix $M_{sX_i}$ of $F(sX_i)$,
$i=1,\cdots,13$, with respect to basis $\{X_1,\cdots,X_{13}\}$ is
\be
&&\hspace{-1cm} \tiny \left[\begin{array}{ccccccccccccc}
1&0&0&0&0&0&0&0&0&0&0&0&0\\0&1&0&0&0&0&0&0&0&0&0&0&0\\0&0&1
&0&0&0&0&0&0&0&0&0&0\\0&0&0&1&0&0&0&0&0&0&0&0&0
\\0&0&0&0&1&0&0&0&0&0&0&0&0\\0&0&0
&0&0&1&0&0&0&0&0&0&0\\0&0&0&0&0&0&1&0&0&0&0&0&0
\\0&0&0&0&0&0&0&1&0&0&0&0&0\\-s&0&0
&0&0&0&0&0&1&0&0&0&0\\0&s&0&0&0&0&0&0&0&1&0&0&0
\\0&0&0&0&0&0&0&0&0&0&1&0&0\\0&0&-
s&0&0&0&0&0&0&0&0&1&0\\0&0&0&0&0&0&0&0&0&0&0&0&1
\end{array}\right],\nonumber \\[3mm]
&&\hspace{-1cm} \tiny \left[\begin{array}{ccccccccccccc}
1&0&0&0&0&0&0&0&0&0&0&0&0
\\0&1&0&0&0&0&0&0&0&0&0&0&0\\0&0&1
&0&0&0&0&0&0&0&0&0&0\\0&0&0&1&0&0&0&0&0&0&0&0&0
\\0&0&0&0&1&0&0&0&0&0&0&0&0\\0&0&0
&0&0&1&0&0&0&0&0&0&0\\0&0&0&0&0&0&1&0&0&0&0&0&0
\\0&0&0&0&0&0&0&1&0&0&0&0&0\\0&-s&0
&0&0&0&0&0&1&0&0&0&0\\-s&0&0&0&0&0&0&0&0&1&0&0&0
\\0&0&-s&0&0&0&0&0&0&0&1&0&0\\0&0&0
&0&0&0&0&0&0&0&0&1&0\\0&0&0&0&0&0&0&0&0&0&0&0&1
\end{array}\right],\nonumber \\[3mm]
&&\hspace{-1cm} \tiny \left[\begin{array}{ccccccccccccc}
1&0&0&0&0&0&0&0&0&0&0&0&0\\0&1&0&0&0&0&0&0&0&0&0&0&0\\0&0&1
&0&0&0&0&0&0&0&0&0&0\\0&0&0&1&0&0&0&0&0&0&0&0&0
\\0&0&0&0&1&0&0&0&0&0&0&0&0\\0&0&0
&0&0&1&0&0&0&0&0&0&0\\0&0&0&0&0&0&1&0&0&0&0&0&0
\\0&0&0&0&0&0&0&1&0&0&0&0&0\\0&0&-
s&0&0&0&0&0&1&0&0&0&0\\0&0&0&0&0&0&0&0&0&1&0&0&0
\\0&s&0&0&0&0&0&0&0&0&1&0&0\\s&0&0
&0&0&0&0&0&0&0&0&1&0\\0&0&0&0&0&0&0&0&0&0&0&0&1
\end{array}\right],
\nonumber \\[3mm]
&&\hspace{-1cm} \tiny \left[\begin{array}{ccccccccccccc}
1&0&0&0&0&0&0&0&0&0&0&0&0\\0&1&0&0&0&0&0&0&0&0&0&0&0\\0&0&1
&0&0&0&0&0&0&0&0&0&0\\0&0&0&1&0&0&0&0&0&0&0&0&0
\\-s&0&0&0&1&0&0&0&0&0&0&0&0\\0&-s
&0&0&0&1&0&0&0&0&0&0&0\\0&0&-s&0&0&0&1&0&0&0&0&0&0
\\0&0&0&-s&0&0&0&1&0&0&0&0&0\\0&0&0
&0&0&0&0&0&1&0&0&0&0\\0&0&0&0&0&0&0&0&0&1&0&0&0
\\0&0&0&0&0&0&0&0&0&0&1&0&0\\0&0&0
&0&0&0&0&0&0&0&0&1&0\\0&0&0&0&0&0&0&0&0&0&0&0&1
\end{array}\right],\nonumber \\[3mm]
&&\hspace{-1cm} \tiny \left[\begin{array}{ccccccccccccc}
1&0&0&0&0&0&0&0&0&0&0&0&0
\\0&1&0&0&0&0&0&0&0&0&0&0&0\\0&0&1
&0&0&0&0&0&0&0&0&0&0\\s&0&0&1&0&0&0&0&0&0&0&0&0
\\0&0&0&0&1&0&0&0&0&0&0&0&0\\0&0&0
&0&0&1&0&0&0&0&0&0&0\\0&0&0&0&0&0&1&0&0&0&0&0&0
\\0&0&0&0&s&0&0&1&0&0&0&0&0\\0&0&0
&0&-s&0&0&0&1&0&0&0&0\\0&0&0&0&0&s&0&0&0&1&0&0&0
\\0&0&0&0&0&0&0&0&0&0&1&0&0\\0&0&0
&0&0&0&-s&0&0&0&0&1&0\\0&0&0&0&0&0&0&0&0&0&0&0&1
\end{array}\right],\nonumber \\[3mm]
&&\hspace{-1cm} \tiny \left[\begin{array}{ccccccccccccc}
1&0&0&0&0&0&0&0&0&0&0&0&0
\\0&1&0&0&0&0&0&0&0&0&0&0&0\\0&0&1
&0&0&0&0&0&0&0&0&0&0\\0&s&0&1&0&0&0&0&0&0&0&0&0
\\0&0&0&0&1&0&0&0&0&0&0&0&0\\0&0&0
&0&0&1&0&0&0&0&0&0&0\\0&0&0&0&0&0&1&0&0&0&0&0&0
\\0&0&0&0&0&s&0&1&0&0&0&0&0\\0&0&0
&0&0&-s&0&0&1&0&0&0&0\\0&0&0&0&-s&0&0&0&0&1&0&0&0
\\0&0&0&0&0&0&-s&0&0&0&1&0&0\\0&0&0
&0&0&0&0&0&0&0&0&1&0\\0&0&0&0&0&0&0&0&0&0&0&0&1
\end{array}\right],
\nonumber \\[3mm]
&&\hspace{-1cm} \tiny \left[\begin{array}{ccccccccccccc}
1&0&0&0&0&0&0&0&0&0&0&0&0\\0&1&0&0&0&0&0&0&0&0&0&0&0\\0&0&1
&0&0&0&0&0&0&0&0&0&0\\0&0&s&1&0&0&0&0&0&0&0&0&0
\\0&0&0&0&1&0&0&0&0&0&0&0&0\\0&0&0
&0&0&1&0&0&0&0&0&0&0\\0&0&0&0&0&0&1&0&0&0&0&0&0
\\0&0&0&0&0&0&s&1&0&0&0&0&0\\0&0&0
&0&0&0&-s&0&1&0&0&0&0\\0&0&0&0&0&0&0&0&0&1&0&0&0
\\0&0&0&0&0&s&0&0&0&0&1&0&0\\0&0&0
&0&s&0&0&0&0&0&0&1&0\\0&0&0&0&0&0&0&0&0&0&0&0&1
\end{array}\right],\nonumber \\[3mm]
&&\hspace{-1cm} \tiny \left[\begin{array}{ccccccccccccc}
1&0&0&0&0&0&0&0&0&0&0&0&0\\0&1&0&0&0&0&0&0&0&0&0&0&0\\0&0&1
&0&0&0&0&0&0&0&0&0&0\\0&0&0&e^s&0&0&0&0&0&0&0&0&0
\\0&0&0&0&e^{-s}&0&0&0&0&0&0&0&0\\0&0&0&0
&0&e^{-s}&0&0&0&0&0&0&0\\0&0&0&0&0&0&e^{-s}&0&0&0&0&0&0
\\0&0&0&0&0&0&0&1&0&0&0&0&0\\0&0&0
&0&0&0&0&0&1&0&0&0&0\\0&0&0&0&0&0&0&0&0&1&0&0&0
\\0&0&0&0&0&0&0&0&0&0&1&0&0\\0&0&0
&0&0&0&0&0&0&0&0&1&0\\0&0&0&0&0&0&0&0&0&0&0&0&1
\end{array}\right],\nonumber \\[3mm]
&&\hspace{-1cm}  \tiny \left[\begin{array}{ccccccccccccc}
e^s&0&0&0&0&0&0&0&0&0&0&0&0
\\0&e^s&0&0&0&0&0&0&0&0&0&0&0\\0&0&e^s&0
&0&0&0&0&0&0&0&0&0\\0&0&0&1&0&0&0&0&0&0&0&0&0
\\0&0&0&0&e^s&0&0&0&0&0&0&0&0\\0&0&0&0
&0&e^s&0&0&0&0&0&0&0\\0&0&0&0&0&0&e^s&0&0&0&0&0&0
\\0&0&0&0&0&0&0&1&0&0&0&0&0\\0&0&0
&0&0&0&0&0&1&0&0&0&0\\0&0&0&0&0&0&0&0&0&1&0&0&0
\\0&0&0&0&0&0&0&0&0&0&1&0&0\\0&0&0
&0&0&0&0&0&0&0&0&1&0\\0&0&0&0&0&0&0&0&0&0&0&0&1
\end{array}\right],
\nonumber \\[3mm]
&&\hspace{-1cm} \tiny \left[\begin{array}{ccccccccccccc} \cos
s&-\sin s&0&0&0&0&0&0&0&0&0&0&0
\\\sin s&\cos s&0&0&0&0&0&0&0&0&0&0&0\\0&0&1&0&0
&0&0&0&0&0&0&0&0\\0&0&0&1&0&0&0&0&0&0&0&0&0
\\0&0&0&0&\cos s&-\sin s&0&0&0&0&0&0&0\\0&0&0&0
&\sin s&\cos s&0&0&0&0&0&0&0\\0&0&0&0&0&0&1&0&0&0&0&0&0
\\0&0&0&0&0&0&0&1&0&0&0&0&0\\0&0&0
&0&0&0&0&0&1&0&0&0&0\\0&0&0&0&0&0&0&0&0&1&0&0&0
\\0&0&0&0&0&0&0&0&0&0&\cos s&\sin s&0\\0&0&0&0&0
&0&0&0&0&0&-\sin s&\cos s&0\\0&0&0&0&0&0&0&0&0&0&0&0&1\end{array}
\right],\nonumber \\[3mm]
&&\hspace{-1cm} \tiny \left[\begin{array}{ccccccccccccc}
1&0&0&0&0&0&0&0&0&0&0&0&0
\\0&\cos s&\sin s&0&0&0&0&0&0&0&0&0&0\\0&-\sin s&\cos s&0
&0&0&0&0&0&0&0&0&0\\0&0&0&1&0&0&0&0&0&0&0&0&0
\\0&0&0&0&1&0&0&0&0&0&0&0&0\\0&0&0
&0&0&\cos s&\sin s&0&0&0&0&0&0\\0&0&0&0&0&-\sin s&\cos
s&0&0&0&0&0&0
\\0&0&0&0&0&0&0&1&0&0&0&0&0\\0&0&0
&0&0&0&0&0&1&0&0&0&0\\0&0&0&0&0&0&0&0&0&\cos s&0&-\sin
s&0\\0&0&0&0&0&0&0&0&0&0&1&0&0\\0&0&0 &0&0&0&0&0&0&\sin s&0&\cos
s&0\\0&0&0&0&0&0&0&0&0&0&0&0&1
\end{array}\right],\nonumber \\[3mm]
&&\hspace{-1cm} \tiny \left[\begin{array}{ccccccccccccc} \cos
s&0&\sin s&0&0&0&0&0&0&0&0&0&0
\\0&1&0&0&0&0&0&0&0&0&0&0&0\\-\sin s&0&\cos s
&0&0&0&0&0&0&0&0&0&0\\0&0&0&1&0&0&0&0&0&0&0&0&0
\\0&0&0&0&\cos s&0&\sin s&0&0&0&0&0&0\\0&0&0&0&0
&1&0&0&0&0&0&0&0\\0&0&0&0&-\sin s&0&\cos s&0&0&0&0&0&0
\\0&0&0&0&0&0&0&1&0&0&0&0&0\\0&0&0
&0&0&0&0&0&1&0&0&0&0\\0&0&0&0&0&0&0&0&0&\cos s&\sin s&0&0
\\0&0&0&0&0&0&0&0&0&-\sin s&\cos s&0&0\\0&0&0&0&0
&0&0&0&0&0&0&1&0\\0&0&0&0&0&0&0&0&0&0&0&0&1
\end{array}\right],\nonumber
 \ee
and $I_{13}$, respectively. Let $X=\sum_{i=1}^{13}a_iX_i$, then
\be && \hspace{-1cm} F(s_7.X_7)\circ F(s_6.X_6)\circ\cdots\circ
F(s_1.X_1)\;:\;X\;\mapsto\;\nonumber\\
&&(s_3a_{12}-s_2a_{10}-s_1a_9+a_1-s_5s_4a_8+s_5a_4-s_4a_5).X_1\nonumber\\
&&+(-s_6s_4a_8+s_3a_{11}-s_2a_9+s_1a_{10}+a_2+s_6a_4-s_4a_6).X_2\nonumber\\
&&+(-s_7s_4a_8-s_3a_9-s_2a_{11}-s_1a_{12}+a3+s_7a_4-s_4a_7).X_3\\
&&+(a_4-s_4a_8).X_4+(-s_6a_{10}-s_5a_9+s_5a_8+a_5+s_7a_{12}).X_5\nonumber\\
&&+(-s_6a_9+s_6a_8+s_5a_{10}+a_6+s_7a_{11}).X_6\nonumber\\
&&+(s_7a_8-s_6a_{11}-s_5a_{12}+a_7-s_7a_9).X_7\nonumber\\
&&+a_8.X_8+a_9.X_9+a_{10}.X_{10}+a_{11}.X_{11}+a_{12}.X_{12}+a_{13}.X_{13}\nonumber
\ee

If $a_8$, and $a_9\neq0$, and $a_8\neq2a_9$, then we can make the
coefficients of $X_1,\cdots,X_7$ vanish. Scaling $X$ if necessary,
we can assume that $a_8=1$. And $X$ is reduced to Case 1.

If $a_8=a_9\neq0$, and $a_{10}\neq0$, then we can make the
coefficients of $X_1\cdots,X_6$ vanish. Scaling $X$ if necessary,
we can assume that $a_8=1$. And $X$ is reduced to Case 2.

If $a_8=a_9\neq0$, $a_{10}=0$, and $a_{12}\neq0$, then we can make
the coefficients of $X_1,\cdots,X_5$, and $X_{10}$ vanish. Scaling
$X$ if necessary, we can assume that $a_8=1$. And $X$ is reduced
to Case 3.

If $a_8=a_9\neq0$, $a_{10}=a_{12}=0$, and $a_{11}\neq0$, then we
can make the coefficients of $X_1,\cdots,X_4$, $X_6$, $X_7$,
$X_{10}$, and $X_{12}$ vanish. Scaling $X$ if necessary, we can
assume that $a_8=1$. And $X$ is reduced to Case 4.

If $a_8=a_9\neq0$, $a_{10}=a_{11}=a_{12}=0$, then we can make the
coefficients of $X_1,\cdots,X_4$, and $X_{10},\cdots,X_{11}$
vanish. Scaling $X$ if necessary, we can assume that $a_8=1$. And
$X$ is reduced to Case 5.

If $a_8=0$, and $a_9\neq0$, then we can make the coefficients of
$X_1,\cdots,X_3$, and $X_5,\cdots,X_8$ vanish. Scaling $X$ if
necessary, we can assume that $a_9=1$. And $X$ is reduced to Case
6.

If $a_8\neq0$, $a_9=0$, and $a_{11}\neq0$, then we can make the
coefficients of $X_1,\cdots,X_3$, and $X_6$, $X_7$, and $X_9$
vanish. And $X$ is reduced to Case 7.

If $a_8\neq0$, $a_9=a_{11}=0$, and $a_{12}\neq0$, then we can make
the coefficients of $X_1,\cdots,X_3$, and $X_5$, $X_7$, $X_9$, and
$X_{11}$ vanish. And $X$ is reduced to Case 8.

If $a_8\neq0$, $a_9=a_{11}=a_{12}=0$, and $a_{10}\neq0$, then we
can make the coefficients of $X_1,\cdots,X_3$, and $X_5$, $X_6$,
$X_9$, and $X_{11}$ vanish. And $X$ is reduced to Case 9.

If $a_8\neq0$, and $a_9=a_{10}=a_{11}=a_{12}=0$, then we can make
the coefficients of $X_1,\cdots,X_3$, and $X_9,\cdots,X_{12}$
vanish. And $X$ is reduced to Case 10.

If $a_8=a_9=0$, $a_4\neq0$, and $a_{10}\neq0$, then we can make
the coefficients of $X_1,\cdots,X_3$, $X_5$, $X_6$, $X_8$, and
$X_9$ vanish. Scaling $X$ if necessary, we can assume that
$a_4=1$. And $X$ is reduced to Case 11.

If $a_4=a_8=a_9=0$, $a_{10}\neq0$, $a_{12}\neq0$, and
$a_6a_{12}+a_7a_{10}\neq a_5a_{11}$, then we can make the
coefficients of $X_1,\cdots,X_6$, $X_8$, and $X_9$ vanish. Scaling
$X$ if necessary, we can assume that $a_7=1$. And $X$ is reduced
to Case 12.

If $a_4=a_8=a_9=0$, $a_{10}\neq0$, $a_{12}\neq0$, and
$a_6a_{12}+a_7a_{10}=a_5a_{11}$, then we can make the coefficients
of $X_1$, $X_2$, and $X_4,\cdots,X_9$ vanish. And $X$ is reduced
to Case 13.

If $a_8=a_9=a_{10}=0$, and $a_4\neq0$, then we can make the
coefficients of $X_1,\cdots,X_3$, $X_5$, and $X_7,\cdots,X_{10}$
vanish. Scaling $X$ if necessary, we can assume that $a_4=1$. And
$X$ is reduced to Case 14.

If $a_8=\cdots=a_{10}=0$, $a_{12}=0$, $a_4\neq0$, and
$a_{11}\neq0$, then we can make the coefficients of
$X_1,\cdots,X_3$, and $X_8,\cdots,X_{10}$ vanish. Scaling $X$ if
necessary, we can assume that $a_4=1$. And $X$ is reduced to Case
15.

If $a_8=\cdots=a_{12}=0$, and $a_4\neq0$, then we can make the
coefficients of $X_1,\cdots,X_3$, and $X_8,\cdots,X_{10}$ vanish.
Scaling $X$ if necessary, we can assume that $a_4=1$. And $X$ is
reduced to Case 16.

If $a_4=a_8=a_9=a_{10}=0$, $a_{12}\neq0$, and $a_5a_{11}\neq
a_6a_{12}$, then we can make the coefficients of $X_1,\cdots,X_5$,
and $X_7,\cdots,X_{10}$ vanish. Scaling $X$ if necessary, we can
assume that $a_6=1$. And $X$ is reduced to Case 17.

If $a_4=a_8=a_9=a_{10}=0$, $a_{12}\neq0$, $a_5a_{11}=a_6a_{12}$,
and $a_6\neq0$, then we can make the coefficients of
$X_1,\cdots,X_5$, and $X_7,\cdots,X_{10}$ vanish. And $X$ is
reduced to Case 18.

If $a_4=a_6=a_8=a_9=a_{10}=0$, $a_{12}\neq0$,
$a_5a_{11}=a_6a_{12}$, and $a_6\neq0$, then we can make the
coefficients of $X_1,\cdots,X_5$, and $X_7,\cdots,X_{10}$ vanish.
And $X$ is reduced to Case 18.

If $a_4=a_6=a_8=\cdots=a_{11}=0$, $a_{12}\neq0$,
$a_5a_{11}=a_6a_{12}$, $a_6\neq0$, and $a_5\neq0$, then we can
make the coefficients of $X_1$, and $X_3,\cdots,X_{11}$ vanish.
And $X$ is reduced to Case 19.

If $a_4=a_8=a_9=a_{10}=a_{12}=0$, $a_5\neq0$, and $a_{11}\neq0$,
then we can make the coefficients of $X_1,\cdots,a_4$,
$X_6,\cdots,X_{10}$, and $a_12$ vanish. Scaling $X$ if necessary,
we can assume that $a_5=1$. And $X$ is reduced to Case 20.

If $a_4=a_8=\cdots=a_{12}=0$, and $a_5\neq0$, then we can make the
coefficients of $X_1$, $a_4$, and $X_8,\cdots,X_{12}$ vanish.
Scaling $X$ if necessary, we can assume that $a_5=1$. And $X$ is
reduced to Case 21.

If $a_4=a_5=a_8=a_9=a_{10}=a_{12}=0$, and $a_{11}\neq0$, then we
can make the coefficients of $X_2,\cdots,X_{10}$ and $a_{12}$
vanish. Scaling $X$ if necessary, we can assume that $a_{11}=1$.
And $X$ is reduced to Case 22.

If $a_4=a_5=a_8=\cdots=a_{12}=0$, and $a_6\neq0$, then we can make
the coefficients of $X_2$, $a_5$, and $X_7,\cdots,X_{12}$ vanish.
Scaling $X$ if necessary, we can assume that $a_6=1$. And $X$ is
reduced to Case 23.

If $a_4=a_5=a_6=0$, and $a_8=\cdots=a_{12}=0$, then we can make
the coefficients of $X_4,\cdots,X_{12}$ vanish. And $X$ is reduced
to Case 24.

If $a_4=a_8=a_9=a_{12}=0$, $a_{10}\neq0$, and $a_7a_{10}\neq
a_5a_{11}$, then we can make the coefficients of $X_1,\cdots,X_4$,
$X_6,\cdots,X_9$, and $a_{12}=0$ vanish. Scaling $X$ if necessary,
we can assume that $a_5=1$. And $X$ is reduced to Case 25.

If $a_4=a_8=a_9=a_{12}=0$, $a_7a_{10}=a_5a_{11}$, and
$a_{10}\neq0$, then we can make the coefficients of
$X_1,\cdots,X_4$, $X_6,\cdots,X_9$, and $a_{12}=0$ vanish. Scaling
$X$ if necessary, we can assume that $a_{10}=1$. And $X$ is
reduced to Case 26. And the Theorem follows. \hfill\ $\Box$

\medskip

According to our optimal system of one-dimensional subalgebras of
the full symmetry algebra $\goth g$, we need only find
group-invariant solutions for 26 one-parameter subgroups generated
by $X$ as Theorem 6. For example, as a direct consequence of
Theorem 3 and case 1 of Theorem 6, we deduce the following
conclusion:
\paragraph{Conclusion 2}
{\em Let $s,c_1,\cdots,c_6\in\RR$ with $c_1$, $c_2$, and
$c_6\neq0$, and $s$ be sufficiently small. Then, if
\begin{eqnarray*}
&& u=U(t,x,y,z),\;v=V(t,x,y,z),\;w=W(t,x,y,z),\\
&& p=P(t,x,y,z),\;q=Q(t,x,y,z),
\end{eqnarray*}
is a solution of the system (\ref{eq:1}), so are the functions
\begin{eqnarray*}
\tilde{u}&=&c_1^sc_2^{-s}\big(\cos (c_3s)\cos (c_5s)-\sin
(c_3s)\sin (c_4s)\sin
(c_5s)\big).U(\tilde{t},\tilde{x},\tilde{y},\tilde{z})\\
&&-c_1^sc_2^{-s}\big(\sin (c_3s)\cos (c_5s)+\cos (c_3s)\sin
(c_4s)\cos
(c_5s)\big).V(\tilde{t},\tilde{x},\tilde{y},\tilde{z})\\
&&+c_1^sc_2^{-s}\cos (c_4s)\sin (c_5s).W(\tilde{t},\tilde{x},\tilde{y},\tilde{z}),\\[2mm]
\tilde{v}&=&c_1^sc_2^{-s}\sin (c_3s)\cos
(c_4s).U(\tilde{t},\tilde{x},\tilde{y},\tilde{z})\\
&&+c_1^sc_2^{-s}\cos (c_3s)\cos (c_4s).V(\tilde{t},\tilde{x},\tilde{y},\tilde{z})\\
&&+c_1^sc_2^{-s}.\sin
(c_4s).W(\tilde{t},\tilde{x},\tilde{y},\tilde{z}),\\[2mm]
\tilde{w}&=&-c_1^sc_2^{-s}\big(\cos (c_3s)\sin (c_5s)+\sin
(c_3s)\sin (c_4s)\cos
(c_5s)\big).U(\tilde{t},\tilde{x},\tilde{y},\tilde{z})\\
&&+c_1^sc_2^{-s}\big(\sin (c_3s)\sin (c_5s)-\cos (c_3s)\sin
(c_4s)\cos
(c_5s)\big).V(\tilde{t},\tilde{x},\tilde{y},\tilde{z})\\
&&+c_1^sc_2^{-s}\cos (c_4s)\cos
(c_5s).W(\tilde{t},\tilde{x},\tilde{y},\tilde{z}),\\[2mm]
\tilde{p}&=&c_6^{-s}.P(\tilde{t},\tilde{x},\tilde{y},\tilde{z}),\\[2mm]
\tilde{q}&=&c_1^{-2s}c_2^{2s}c_6^{-s}.Q(\tilde{t},\tilde{x},\tilde{y},\tilde{z}),
\end{eqnarray*}
where $\tilde{t}=c_1^s.t$, and
\begin{eqnarray*}
\tilde{x}&=& c_2^s\big(\sin (c_5s)\sin (c_4s)\sin (c_3s)+\cos
(c_5s)\cos
(c_3s)\big).x\\
&&+c_2^s\big(\sin (c_3s)\cos (c_5s)-\cos (c_3s)\sin (c_4s)\sin
(c_5s) \big).y\\
&&-c_2^s\cos (c_4s)\sin (c_5s).z,\\[2mm]
\tilde{y}&=& -c_2^s\sin (c_3s)\cos (c_4s).x+(c_2s)\cos
(c_3s)\cos (c_4s).y-c_2^s\sin (c_4s).z,\\[2mm]
\tilde{z}&=& c_2^s\big(\cos (c_3s)\sin (c_5s)-\sin (c_3s)\sin
(c_4s)\cos
(c_5s)\big).x\\
&&+c_2^s\big(\cos (c_3s)\sin (c_4s)\cos (c_5s)+\sin (c_3s)\sin
(c_5s)\big).y\\
&&+c_2^s\cos (c_4s)\cos (c_5s).z.
\end{eqnarray*}}
The construction of the group-invariant solutions for each of the
one-dimensional subgroups in the optimal system proceeds in the
same fashion.
\section{Acknowledgements}
I am grateful to professor Ian M. Anderson for his consideration
and his mathematical career to the development the
"DifferentialGeometry" package of MAPLE 11. It is an applicable
tool to extremely complex and tedious geometrical computations.

\end{document}